\documentclass[a4,invmat]{amsart}
\usepackage{graphicx}
\usepackage{amssymb}
\usepackage{enumerate}
\newtheorem{thm}{Theorem}
\newtheorem{cor}{Corollary}
\newtheorem{lem}{Lemma}
\newtheorem{prop}{Proposition}
\newtheorem{rem}{Remark}

\newcommand{\norm}[1]{\left\Vert#1\right\Vert}

\newcommand{\set}[1]{\left\{#1\right\}}
\newcommand{\Real}{\mathbb R}

\newcommand\Dol{\mathcal{D}}
\newcommand\Fol{\mathcal{F}}

\newcommand\TT{\mathbb{T}}

\newcommand\beqas{\begin{eqnarray*}}
\newcommand\eeqas{\end{eqnarray*}}
\newcommand{\beq}{\begin{equation}}
\newcommand{\eeq}{\end{equation}}

\newcommand{\normsup}[1]{\| #1 \|_\infty}
\newcommand{\normLone}[1]{\| #1 \|_1}

\newcommand\OCInt[1]{(#1]}

\newcommand\One{\mathbf{1}}
\newcommand\Zero{\mathbf{0}}
\newcommand\LL{\mathbf{L}}
\newcommand\UU{\mathbf{U}}

\usepackage{latexsym}
\usepackage{graphics}
\usepackage{amsmath}
\usepackage{xspace}
\usepackage{amssymb}

\newcommand {\be}[1]{\begin{equation}\label{#1}}
\newcommand {\ee}{\end{equation}}
\newcommand {\bea}{\begin{eqnarray}}
\newcommand {\eea}{\end{eqnarray}}
\newcommand{\Pois}{\ensuremath{\operatorname{Pois}}\xspace}

\newcommand{\pr}{\mathbb{P}}
\newcommand{\E}{\mathbb{E}}

\begin{document}
\title{Dynamics of exponential linear map in functional space}
\author{David Gamarnik}
\author{Tomasz Nowicki}
\author{Grzegorz \'Swirszcz}
 \address{IBM Watson Research Center, Yorktown Heights NY 10598, USA.}
 \email{gamarnik@watson.ibm.com}
 \email{tnowicki@watson.ibm.com}
 \address{Institute of Mathematics, University of Warsaw, 02--097 Warsaw, Banacha 2, Poland.}
 \address{IBM Watson Research Center, Yorktown Heights NY 10598, USA.}
 \email{swirszcz@mimuw.edu.pl}
 \email{swirszcz@us.ibm.pl}
 \thanks{GS was partially supported by Polish KBN Grant 2PO3A 01022}

\maketitle
\begin{abstract}
 We consider the question of existence of a unique invariant
 probability  distribution which satisfies some evolutionary
 property. The  problem arises from the random graph theory but to
 answer it we treat it as a dynamical system in the functional
 space, where we  look for a global attractor. We consider the
 following bifurcation problem: Given a probability measure $\mu$,
 which corresponds to the  weight distribution of a link of a
 random graph we form a positive linear operator $\Phi$
 (convolution) on distribution functions  and then we analyze a
 family of its exponents with a parameter $\lambda$ which
 corresponds to connectivity of a sparse random graph. We prove
 that for every measure $\mu$ (\emph{i.e.}, convolution $\Phi$)
 and every $\lambda< e$ there exists a unique globally attracting
 fixed point of the operator, which yields the existence and
 uniqueness of the limit probability distribution  on the random
 graph. This estimate  was established earlier \cite{KarpSipser}
 for deterministic weight distributions (Dirac measures $\mu$) and
 is known as $e$-cutoff phenomena, as for such distributions and
 $\lambda>e$ there is no fixed point attractor. We thus establish
 this phenomenon  in a much more general sense.
\end{abstract}

\section{Introduction}\label{sec:Intr}
A dynamical system is a model of  time  evolution. If the
asymptotic behavior of the system is independent on the initial
conditions then we can say that the system forgets about its past,
or that it is impossible to reconstruct the past knowing the far
future. A simplest such situation arises when the system has
\emph{a fixed point which is a global attractor}, in other words
that wherever we started our trajectory we land in the same spot.
Our paper was motivated by studying this approach to some aspects
of the theory of random graphs, which we will explain in some
details after the definitions. The importance of the uniqueness of
a fixed point of the dynamical system is related to the effect of
decay of correlation in the underlying random graph. Specifically,
if the fixed point is unique then structure in one part of the
graph is asymptotically independent from such a structure in other
parts of the graph. The connection between uniqueness of a fixed
point and correlation was formally established by the authors in
\cite{gamarnikMaxWeightIndSet}. The concept of correlation decay
comes up frequently in statistical physics. In a particular
context of Glauber dynamics on spin glasses on trees see
\cite{Martin}, \cite{MartinellySinclairWeitzMixing} (also
\cite{BergerKenyonMosselPeres}, \cite{BrightwellWinkler}, and
related problem of information flow on trees
\cite{EMosselSurvey}). In dynamical systems
 it is often connected with the existence of a
unique invariant measure (or an attracting  fixed point of a
Perron Frobenius operator) \cite{Baladi:DecayBook}.

The intention of this paper is to link different fields: dynamics,
probability, graphs and analysis. We provide therefore detailed
proofs, to make the results  accessible for readers with different
background.

We consider a class of nonlinear operators in the space of
distributions which arises in the study of maximum weight
matchings in sparse random graphs.  The study of this object leads
to the problem of existence and uniqueness of an invariant
distribution for an iterative process and whether the iterations
of any distribution converges to this invariant distribution. Our
methods rely on the understanding of the dynamics of this
iterative process. The phase space is however a functional space
and we have a rare opportunity to study a specific non-linear
system with non-trivial behavior. The nonlinear operator is a
composition of an exponential map with a positive linear operator.

\subsection*{Definitions and main result}\label{sub:Definitions}
Let $\Fol$ be a family of functions on the segment $I\subset
\Real$ with values in $[0,1]$, $\Fol\ni F:I\to\Real$. For a
positive linear operator (endomorphism) $\Psi$ on $\Fol$ we
define:
\[
\TT:\Fol\to\Fol\qquad\text{by}\qquad\TT(F)(x)=\exp(-\Psi(F)(x))\,.
\]
In the case which interests us at most as it has an application to
the random graphs the linear operator $\Psi$ is a  product of the
parameter $\lambda>0$ and a convolution $\Phi$ with respect to a
given probability measure $\mu$.
 We restrict the domain of $\Phi$  to nondecreasing functions of an
interval, which we fix here to be $[0,1]$. Specifically: $\Dol\ni
F(z):[0,1]\to[0,1] $ and  $F$ is not decreasing. Given a (Borel)
probability measure $\mu$ on $[0,1]$ we define a linear operator
on $\mu$ integrable functions (Lebesgue integral) by:
\[
 \Phi_\mu(F)=\int_x^1F(z-x)\,d\mu(z)=\int_0^{1-x}F(z)\,d\mu(z+x)\,,
 \]
 and with it, for any $\lambda>0$, an exponential map by:
\[
  F\mapsto \TT_\lambda(F)(x)=\TT F(x)=\exp\left(-\lambda
  \Phi_\mu(F)(x)\right)\,.
\]
The main result of this paper is an extension of the
\emph{$e$-cutoff phenomenon}, which was established earlier by
Karp and Sipser \cite{KarpSipser} only for deterministic
distributions. We state this result in three versions:
\begin{itemize}
    \item Theorem~\ref{thm:contraction for T^2}, which carries the burden of the proof. For any positive
        linear $\Psi$ : If $\norm{\Psi}<e$, then $\TT^2$ is a contraction.
    \item Theorem~\ref{thm:TheECut} is convolution specific. For any probability measure $\mu$,
        if $\lambda<e$, then $\TT$ has a fixed point  which is a global attractor.
    \item Theorem~\ref{thm:ECutRandomGraph} restates the result in
    the probabilistic setting.
\end{itemize}
\subsection*{Applications to maximum weight matching in sparse
random graphs}\label{sub:RandomGraph} Before we prove our main
results we  describe in more details and in the the probabilistic
setting the connection between the fixed point properties of the
dynamical systems considered in Sections \ref{sec:ExpoLinea} and
\ref{sec:Distributions} and the theory of random graphs.

The following is a standard model of a sparse random graph on $n$
nodes with average degree (connectivity) $\lambda$
\cite{BollobasBook}, \cite{JansonBook}. Often this model is also
called Erdos-Renyi graph. Given a collection of $n$ nodes
$1,2,\ldots,n$, an edge (link) $(i,j), 1\leq i\neq j\leq n$ is
selected to belong to the  graph with probability $\lambda/n$,
independently for all $n(n-1)/2$ pairs $i,j$. The collection of
selected edges is denoted by $E$. The selected edges are equipped
with randomly generated non-negative weights $W_{i,j}$,
distributed according to a common distribution function
$\pr(W_{i,j}\leq x)\equiv \mu(x), x\geq 0$. A matching is any
collection of edges in $E$ which do not share a node. That is
$M\subset E$ is a matching if for every $(i_1,j_1),(i_2,j_2)\in E$
the nodes $i_1,i_2,j_1,j_2$ are distinct. The weight of a matching
$M$ is the sum $\sum_{(i,j)\in M}W_{i,j}$. We let
$M_{\mu}(n,\lambda)=\max_M\sum_{i,j}W_{i,j}$ denote the maximum
weight of a matching. Note that $M_{\mu}(n,\lambda)$ is a random
variable which only depends on $n,\lambda$ and the distribution
function $\mu$. The main question of interest is establishing the
existence and computing the limit
\begin{equation}\label{eq:MatchingLimit}
\lim_{n\to\infty}{\mathbb{E}[M_{\mu}(n,\lambda)]\over n}\,,
\end{equation}
 where $\mathbb{E}[\cdot]$ denotes the expectation operator.
This problem was solved for the case of the deterministic weights
($W_{i,j}=1$ with probability one (w.p.1)) by Karp and Sipser
\cite{KarpSipser} using a simple combinatorial argument, that we
reproduce for completeness in the third part of
Theorem~\ref{thm:examples}.
The
threshold $\lambda=e$ also corresponds to some phase
transition property in the underlying random graph. This phase
transition  was thereafter called $e$-cutoff phenomena.

The Karp-Sipser method however does not apply  to the case of
non-deterministic weights and the authors
\cite{gamarnikMaxWeightIndSet} solved the problem of computing
limit \eqref{eq:MatchingLimit} using the completely different
Local Weak Convergence (LWC) method, developed earlier by Aldous
\cite{Aldous:assignment92}, \cite{Aldous:assignment00}, Aldous and
Steele\cite{AldousSteele:survey}. The existence of the limit
\eqref{eq:MatchingLimit} was first established by the first author
\cite{gamarnik_LSAT} using a non-constructive version of LWC
method and only in \cite{gamarnikMaxWeightIndSet} we were able to
compute the limit at least for some non-deterministic
distribution. The method is heavily based on solving for fixed
point solutions of certain distributional equations. We give here
only a quick description of the main result in
\cite{gamarnikMaxWeightIndSet} regarding matching and refer the
reader to the paper for further details.

Let $F=F(x)$ be a distribution function corresponding to some
non-negative random variable, and let $K$ be a random variable
distributed according to the Poisson distribution with parameter
$\lambda$, denoted $\Pois(\lambda)$. That is
$\pr(K=k)=(\lambda^k/k!)e^{-\lambda}$. Consider a random variable
$X=\max_{i\leq K}(W_i-X_i)$, where $X_1,\ldots,X_K$ are
distributed according to $F$, independently and $W_1,\ldots,W_K$
are distributed according to $\mu=\mu(x)$ independently. When
$K=0$, $X$ is assumed  $0$ by convention. Let $\tilde F$ denote
the distribution function of $X$. This defines an operator
$F\mapsto \tilde F$ on the space of distribution functions,
indexed by $\lambda$ and the distribution function $\mu$. We claim
that this operator is in fact $F\mapsto\TT_\lambda(F)$ defined in
Section \ref{sec:Distributions}. Indeed:
    \beqas
    \pr(X\leq x)
    & = & \sum_{k=0}^{\infty}{\lambda^k\over k!}e^{-\lambda}\left(\pr(W_1-X_1\leq
    x)\right)^k
    \\
    & = & \sum_{k=0}^{\infty}{\lambda^k\over k!}e^{-\lambda}
        \left(\mu(x)+\int_x^{\infty}(1-F(z-x))d\mu(z)\right)^k
    \\
    & = & e^{-\lambda}\sum_{k=0}^{\infty}{\lambda^k\over k!}
        \left(\mu(\infty)-\int_x^{\infty}F(z-x)d\mu(z)\right)^k
    \\
    & = & \exp\left(-\lambda\int_x^{\infty}F(z-x)d\mu(z)\right)\,,
    \eeqas
where in the last equality we use the fact that $\mu$ is the
distribution function, and therefore $\mu(\infty)=1$. We see that
the distribution $\tilde F$ of $X$ is indeed given by
$\TT_\lambda(F)$. The following theorem was established in
\cite{gamarnikMaxWeightIndSet} (Theorem 2, Equation (9)).

\begin{thm}\label{thm:MainRandomGraph}
Given an atom-free distribution function $\mu$, suppose the
operator $\TT_\lambda^2$ has the unique fixed point solution
$F^*=\TT_\lambda^2(F^*)$. Then the limit \eqref{eq:MatchingLimit}
is equal to:
 \begin{equation}\label{eq:expression}
    {1\over 2}\E\left[\sum_{i\leq
    K}W_i\,\chi\{W_i-X_i=\max_{j\leq K}(W_j-X_j)>0\}\right]\,,
 \end{equation}
where $K$ is distributed as $\Pois(\lambda)$, $W_1,\ldots,W_K$ are
distributed according to $\mu$, and $X_1,\ldots,X_K$ are
distributed according to $F^*$. ($\chi A$ is the indicator
function of the set $A$).
\end{thm}
The expression above can easily be transformed into an expression
involving integrals and distribution functions. Thus the theorem
states that whenever the fixed point $F^*$ of $\TT_\lambda^2$ is
unique, the maximum weight matching can be computed by computing
the expectation above with respect to measures $\mu$, $F^*$ and
$\Pois(\lambda)$. The theorem then justifies the search for
measures $\mu$ and parameters $\lambda$ for which the
corresponding operator $\TT_\lambda^2$ has the unique fixed point.
Later on Theorem \ref{thm:TheECut} shows that for every
distribution $\mu$ and every $\lambda<e$ the operator
$\TT_\lambda$ does have a unique fixed point which is a global
attractor and, as a result, $\TT_\lambda^2$ has a unique fixed
point. We then solve the problem of finding
(\ref{eq:MatchingLimit}) whenever $\lambda<e$. Hence the
$e$-cutoff Theorem \ref{thm:TheECut} in a context of random graphs
is:
\begin{thm}[\textbf{$e$-cutoff probabilistic version}]\label{thm:ECutRandomGraph}
Given a sparse random graph with connectivity $\lambda<e$ and
given an atom-free distribution function $\mu$, the limit
\eqref{eq:MatchingLimit} is equal to \eqref{eq:expression}, where
$F^*$ is the unique fixed point of $\TT_\lambda$.
\end{thm}

\subsection*{Outline of the paper}\label{sub:Outline}
\begin{itemize}
    \item In next Section~\ref{sec:ExpoLinea} we
        prove that the exponential linear map defines a dynamical system
        (Proposition~\ref{prop:GeneralTwelldefined}) on real functions and prove its basic
        properties: monotonicity, continuity and differentiability in
        natural norms.
    \item We establish the existence of two specific limit functions
         and in Theorem~\ref{thm:Dichotomy} we
        prove (using only monotonicity) that the system has a fixed  point
        global attractor if and only if they are equal. When the operator is continuous
        in the $L^1$ norm, the limit functions form a periodic cycle and
        Theorem~\ref{thm:JawsInL1} specifies in such a case that the existence of
        a fixed point global attractor is equivalent to the
        uniqueness of the fixed point for $\TT^2$.
    \item We prove that if the sup norm of the linear operator is smaller than $e$, the second iterate
        of the exponential map is a contraction and the main
        result (Theorem~\ref{thm:contraction for T^2})  follows.
    \item In Section~\ref{sec:Distributions} we deal with specific linear part, the convolutions.
        We restrict the phase space to a subset of  non decreasing
        functions (in fact with range in $[0,1]$)
        and prove that in this  case the map defines the dynamical system on this set of
        distribution functions.
    \item We prove a criterion (Theorem~\ref{thm:EndPointCriterion}) for the existence of the
        fixed point global attractor which is specific to the restricted system.
    \item In last Section~\ref{sec:Examples} we present examples of the map for
        particular measures $\mu$. In cases of Lebesgue measure (uniform
        distribution) and exponential distribution there are  fixed
        points which are a global attractors for every $\lambda>0$. For
        completeness we also include  the known case of Dirac measure
        where  the fixed point is a global attractor if and only if $0<\lambda\le e$.
\end{itemize}
\section{The exponential-linear dynamics for positive linear operator $\Psi$}\label{sec:ExpoLinea}
    In order to lighten the notation we will write $\TT F$ and $\Psi F$ for $\TT(F)$ and $\Psi(F)$.
\begin{prop}[\textbf{$\TT$ defines a dynamical system on real functions}]\label{prop:GeneralTwelldefined}
 For every $F\in\Fol$ we have $\TT F\in \Fol$
\end{prop}
\begin{proof}
By positivity of $\Psi$, if  $F\ge 0$ then $\Psi F\ge 0$. Hence
$\exp(-\Psi F)\in[0,1]$.
\end{proof}
\begin{rem}
The definition of $\TT$ can be extended by linearity of $\Psi$ to
any function $F:I\to\Real$ such  that  $a F(x)+b$ lies in the
domain of $\Psi$ for some $a,b\in\Real$, for example to bounded
functions. We have for any  such $F$:
\[
\TT F\ge 0\qquad\text{ and }\qquad \TT(\TT F)\le 1\,,
\]
so our assumption on the range of $F\in\Fol$ is not very
restrictive. \qed
\end{rem}
 Let $\TT^n$ denote the $n$-th iteration of $\TT$ given by $\TT^0(F)(x)=F(x)$ and
$\TT^{n+1}(F)=\TT(\TT^n(F))$.
\subsection*{Monotonicity properties}\label{sub:Monotonicity}
\begin{lem}\label{lem:Monotonicity}
The map $\TT$ is non increasing, the map $\TT^2$ is non
decreasing.\qed
\end{lem}
\begin{proof}
If $F, G\in\Fol$, $F\le G$ then $G-F\in\Fol$ and $0\le \TT(G-F)\le
1$. By linearity of $\Psi$ we have $\TT(G)=\TT
(F)\cdot\TT(G-F)\le\TT(F)$. For $\TT^2$ we apply the previous
argument twice.
\end{proof}
Define:
\begin{equation}\label{eq:ZeroOne}
\Zero(x)=0\quad\text{and}\quad\One(x)=1\quad\text{for all}\quad
x\in I\,,
\end{equation}
we have $\TT(\Zero)=\One$. Denote:
\[
\Zero^n=\TT^n(\Zero)\quad\text{and}\quad\One^n=\TT^n(\One)\,,
\]
we have clearly $\Zero^{n+1}=\One^n$, which symbolically defines
$\One^{-1}=\Zero$.
\begin{lem}\label{lem:Jaws}
For every $F\in\Fol$ and every $n\ge 0$ we have:
 \begin{equation*}
 \begin{array}{rcl}
  \Zero^{2n}=\One^{2n-1}\le&\TT^{2n}F&\le \One^{2n}\\
\text{and}\qquad  \One^{2n}\ge &\TT^{2n+1}F&\ge
\One^{2n+1}=\Zero^{2n+2} \,,
 \end{array}
 \end{equation*}
 and in  particular:
\[
\Zero=\One^{-1}\le
\dots\le\One^{2n-1}\le\One^{2n+1}\le\dots\le\One^{2n+2}\le\One^{2n}\le
\dots\le \One^0=\One\,.
\]
\end{lem}
\begin{proof}
By definition $\Zero\le F\le\One$. Hence by
Lemma~\ref{lem:Monotonicity}, $\One=\TT\Zero\ge \TT F\ge\TT
\One=\TT^2\Zero$ and the inequalities follow by induction.
\end{proof}
\begin{cor}\label{cor:UL}
There are point-wise, monotone limits
\[
\LL=\lim \Zero^{2n}\le\lim\One^{2n}=\UU\,.
\]
For every $F\in\Fol$ we have
\[
\LL\le \liminf \TT^n F\le \limsup\TT^n F\le \UU\,.
\]\qed
\end{cor}
\begin{thm}[\textbf{Main criterion for uniqueness of the attractor}]\label{thm:Dichotomy}
The exponential linear dynamical system has  a fixed point which
is a global attractor if and only if the limit functions $\LL$ and
$\UU$ are equal.
\end{thm}
\begin{proof}
If $\LL=\UU$ then by Corollary~\ref{cor:UL} every $\TT^nF$
converges point-wise to a common limit. If $\LL\not=\UU$ then
$\TT^n\One$ (and any function contained between two odd or two
even iterates of $\One$) do not converge to a limit as it has two
distinct accumulation points $\LL$ and $\UU$.
\end{proof}
\begin{rem}
If $\LL\not=\UU$ in  Theorem~\ref{thm:Dichotomy} then the  a
global attractor (which is not a fixed point and may not be
minimal) is contained in the set $\set{F\in\Dol:\LL\le F\le \UU}$
and contains both $\LL$ and $\UU$. It seems that both inclusions
are proper. It still may happen that some functions from between
$\LL$ and $\UU$ converge  to a fixed point.\qed
\end{rem}
\subsection*{Continuity and differentiability}\label{sub:ContinuityAndDifferentiability}\ \\
From now on we assume that the linear operator is continuous in
either the $\sup$ norm $\normsup{\cdot}$ or in the $L^1(dx)$ norm
$\normLone{\cdot}$ on $I$. In the second case, when the segment
$I$ is infinite we assume that every $F\in\Fol$ has a bounded
integral; moreover all equalities of the functions are meant in
the norm sense, \emph{i.e.,} $F=G$ means $\normLone{G-F}=0$.
\begin{lem}\label{lem:Tcontinuous}
If $\Psi$ is continuous (or in other words when its norm is
bounded) then $\TT$ is continuous.
\end{lem}
\begin{proof}
 It follows from the composition rule.
\end{proof}
\begin{rem}\label{rem:Tcontinuous}
 If $\TT$ is continuous in some norm and the sequence $\One^{2n}$ converges in the same
 norm, then $\One^{2n+1}$ converges and the limits are $\UU$ and
 $\LL$ respectively. In such a case:
 \[
\TT\LL=\UU\quad\text{and}\quad\TT\UU=\LL\,.
\]
In the $L^1$ norm the existence of the limit is assured by the
Lebesgue Convergence Theorems (either monotone or majorized).
\end{rem}
\begin{cor}[\textbf{The $L^1$ norm}]\label{thm:JawsInL1}
Suppose that $\TT \UU=\LL$ and $\TT\LL=\UU$. Then $\TT$ has a
fixed point which is a global attractor if and only if $\TT^2$ has
a  unique fixed point.
\end{cor}
\begin{proof}
Both $\LL$ and $\UU$ are fixed points of $\TT^2$ and the result
follows from Theorem~\ref{thm:Dichotomy}.
\end{proof}
\begin{lem}\label{lem:Tdiff}
If $\Psi$ is continuous then the operator $\TT$ is differentiable
with respect to $F$ and its derivative is:
\[
D\TT (F)(H)(x)=\TT F(x)\cdot\Psi H(x)\,.
\]
\end{lem}
\begin{proof}
The linear operator $\Psi$ is continuous and hence differentiable.
The formula is an application of the Chain Rule.
\end{proof}
\begin{cor}\label{lem:normT(F+H)-T(F)}
The derivative of $\TT$ is uniformly bounded for $F\in\Fol$ by the
norm of $\Psi$ and:
\[
\norm{\TT(F+H)-\TT(F)}\le \norm{\Psi}\norm{H}\,.
\]
\end{cor}
\begin{proof}
We have $|\TT F(x)\cdot\Psi H(x)|\le \sup|\TT F|\cdot|\Psi
H(x)|\le |\Psi H(x)|$. The formula follows from the Mean Value
Theorem.
\end{proof}
\begin{cor}\label{cor:contraction for T}
If $\norm{\Psi}<1$ then $\TT$ has a unique fixed point which is a
global attractor.
\end{cor}
\begin{proof}
In this case $\TT$ is a contraction, hence $\UU$ and $\LL$ cannot
stay away and we use Theorem~\ref{thm:Dichotomy}. (Note that we do
not need to prove that $\Fol$ is complete.)
\end{proof}
\begin{thm}\label{thm:contraction for T^2}
If $\normsup{\Psi}< e$ then $\TT$ has a unique fixed point which
is a global attractor.
\end{thm}
\begin{proof}
It is enough to show that $\TT^2$ is a contraction in
$\normsup{\cdot}$ hence then again $\UU$ and $\LL$ cannot be
different. We use the Chain Rule to calculate the derivative of
the second iterate:
\[
D(\TT^2)(F)(H)=\TT^2F\cdot\Psi(\TT F\cdot\Psi H)\,,
\]
Since for $G>0$ and any $k$ we have
$\Psi(G\cdot(\normsup{k}-k))\ge 0$, then $|\Psi(G\cdot k)|\le|\Psi
G|\cdot \normsup{k}$ and therefore:
\[
|D\TT^2(F)(H)|\le \TT^2F\cdot\Psi(\TT F)\cdot\normsup{\Psi H}\le
\frac{\normsup{\Psi}}{e}\cdot\normsup{H}\,,
\]
where we used the fact that $xe^{-x}\le 1/e$ for $0\le x=\Psi(\TT
F)$. Now from the Mean Value Theorem we see that
$\normsup{\TT^2(F+H)-\TT^2(F)}<\kappa\normsup{H}$, where
$\kappa=\normsup{\Psi}/e<1$.
\end{proof}
\begin{rem}\label{rem:GeneralCutoff}
The bound is tight, for $\lambda>e$ the map $\TT$ may have
different type of global attractors, see Theorem~\ref{thm:TheECut}
and  Theorem~\ref{thm:examples} (and
Remark~\ref{rem:DeltaOtherAttractors} in its proof).\qed
\end{rem}
\section{The linear part of the operator $\TT$ is a convolution on the set of distribution functions.}
\label{sec:Distributions} In this section $\Psi=\lambda\Phi$,
where $\Phi$ is a convolution of  nondecreasing functions of an
interval $[0,1]$ with respect to a given probability measure
$\mu$. For $\Dol\ni F(z):[0,1]\to[0,1] $ and $F$ is not
decreasing, $F$ can be extended by zero to the left of $0$ and by
$1$ to the right of $1$. We do not assume that $F(0)=0$, in other
words the measure defined by $\int dF$ may have an atom at $0$. As
it is of no consequence to our result we do not resolve the
continuity issue at jump points. The two particular functions
$\Zero, \One\in\Dol$. Note as curiosity that $\int_0^x
d\mu(z)\in\Dol$.
\begin{rem}\label{rem:AtomAtOne}
 We assume that $\mu$ has no atom at $1$,
 otherwise by rescaling the interval of arguments
 we can push the support of $\mu$ inside $[0,1]$. However if
 $\mu$ has no atom at the right-most point of its support, again
 by rescaling we may assume that the point $1$ belongs to the support
 of the measure $\mu$. In both cases $\int_1^1d\mu=0$, and hence
 we may add the condition $F(1)=1$, which will be preserved by
 $\TT$.\qed
 \end{rem}
\begin{prop}[\textbf{$\TT$ defines a dynamical system on $\Dol$}]\label{prop:TWellDefined}
For $F\in\Dol$ we have  $\TT F\in\Dol$.
\end{prop}
\begin{proof} We know  that if $\Phi$ is positive $\TT F(x)\in[0,1]$.
We have to check that $\TT F(x)$ is
non decreasing and (if we apply the convention of
Remark~\ref{rem:AtomAtOne}) that $\TT F(1)=1$.
\begin{enumerate}
 \item $\Phi$ is a positive linear operator. Clearly if $ F\ge 0$,
  then $\Phi F\ge 0$. Hence for $\lambda>0$ also $0\le \TT F\le 1$.
 \item If $y-x\ge 0$, then by assumption $F(y)- F(x)\ge 0$ and
  $F(z-y)-F(z-x)\le 0$, therefore:
  \beqas
  \Phi F(y)-\Phi F(x)=\int_y^1 (F(z-y)-F(z-x))\,d\mu(z)-\int_x^y
  F(z-x)\,d\mu(z)\le 0\,.
  \eeqas
Hence:
 \[
 \TT F(y)-\TT F(x)=\TT F(y)\cdot\left(1-\exp\left(-\lambda(\Phi F(x)-\Phi F (y))\right)\,\right)\ge
 0\,.
 \]
 \item By  Remark~\ref{rem:AtomAtOne}: $\Phi F(1)=0$, thus $\TT F(1)=\exp(-\lambda \cdot 0)=1$.
\end{enumerate}
\end{proof}
\begin{lem}[\textbf{$\Phi$ is self adjoint}]\label{lem:SelfAdjoint}
The expected value $\E_\mu$ with respect to the measure $\mu$ of
the convolution $F\star G$ with respect to the Lebesgue measure is
equal to:
\[
\int_0^1 F\cdot\Phi G\,dx=\E_\mu [F\star G]=\int_0^1\Phi F\cdot
G\,dx\,.
\]
\end{lem}
\begin{proof}
We use Fubini Theorem and the change of variables $x=z-w$, all
three variables being in $[0,1]$:
 \beqas
&&\int_0^1 F(x)\cdot\Phi G(x)\,dx=
   \iint\limits_{z\ge x} F(x)G(z-x)\,d\mu(z)\,dx\\
&&\qquad=
   \iint\limits_{w\le z}F(z-w)G(w)\,dw\,d\mu(z)
   =\int_0^1\Phi F(w)\cdot G(w)\,d\mu(z)\,dw\,.
\eeqas In fact we could skip the limits as $F(x),G(x)$ are zero
for $x<0$ and the support of $\mu$ is in $[0,1]$.
\end{proof}
\begin{prop}\label{prop:Continuity of Phi}\ \\
The linear operator $\Phi$ is continuous (and hence
differentiable) in the norms $\normsup{\cdot}$ and
       $\normLone{\cdot}$, and its operator norm in both cases does not
       exceed 1.
\end{prop}
\begin{proof}
We see that $\sup{|\Phi(H)|}\le\Phi(\sup|H|)\le \sup|H|$, where we
used that for any  function $K(x)\le  k$ with $k>0$, we have
$\Phi(K)(x)\le k\cdot\mu([x,1])\le k$.  Using this inequality with
$K(x)=\One(x)$ and Lemma~\ref{lem:SelfAdjoint}  we have:
        \[
        \int_0^1|\Phi H|\,dx\le\int_0^1 \One\cdot\Phi|H|\,dx\le\int_0^1 |H|\cdot\Phi\One\,dx
        \le\int_0^1|H|\,dx\,.
        \]
\end{proof}
\begin{rem}\label{rem:Norm=1}
 In fact $\normsup{\Phi}=1$ as can be checked by $H={\rm const}$.
Also  $\normLone{\Phi}=1$ if $\mu$ has no atom at
 one (as then, taking $H_n$ constant on $\OCInt{1-1/n,1}$ leads to the
 needed estimate). If $\mu$ has an atom at one with weight $p$ then
 $\normLone{\Phi}\le 1-p$.
 \qed
\end{rem}
\begin{cor}
 Continuity and differentiability of $\TT$ follow from
 Proposition~\ref{lem:Tcontinuous} and Lemma~\ref{lem:Tdiff}.
\end{cor}
\begin{thm}[\textbf{The
$e$-cutoff}]\label{thm:TheECut}
 For any
$0<\lambda< e$ and any probability measure $\mu$ the map
$\TT_\lambda(\cdot)=\exp(-\lambda\Phi_\mu(\cdot))$ has a unique
fixed point which is a global attractor. On the other hand there
exists $\mu$ such that for $\lambda>e$ the map $\TT$ has no fixed
point, which is a global attractor.
\end{thm}
\begin{proof}
An example  where there is no fixed point global attractor for
$\lambda>e$  is presented in Theorem~\ref{thm:examples} (the case
of Dirac measure), see also Remark~\ref{rem:DeltaOtherAttractors}.
For $\lambda<e$ the result follows from
Theorem~\ref{thm:contraction for T^2} with $\Psi=\lambda\Phi$, as
$\norm{\Psi}=\lambda\norm{\Phi}\le \lambda$.
\end{proof}
Now we present a technical condition which may help to decide on
the existence of a globally attracting fixed point.
\begin{thm}\label{thm:EndPointCriterion}
If there exists an  $N$ such that
\[
\One^{2N+1}(0)>\frac{1}{e}\,,
\]
then $\TT$ has a fixed point which is a global attractor in the
norm $\normLone{\cdot}$.
\end{thm}
Note that we consider  an odd iterate of $\One$ or equivalently an
even iterate of $\Zero$.
\begin{proof}
From the  assumption it follows that $\LL(x)\ge
\LL(0)\ge\One^{2n+1}(0)>1/e$. The function $x\ln(x)$ is increasing
for $x>1/e$ hence, as $1/e\le\LL\le \UU$ we have $\LL\ln\LL\le
\UU\ln\UU$. On the other hand $\ln\LL=-\lambda\Phi\UU$ and
$\ln\UU=-\lambda\Phi\LL$:
\[
0\le \int_0^1(\UU\ln\UU-\LL\ln\LL)\,dx=
-\lambda\int_0^1(\UU\cdot\Phi\LL-\LL\cdot\Phi\UU)\,dx=0\,,
\]
by Lemma~\ref{lem:SelfAdjoint}. Hence $\normLone{\UU-\LL}=0$ and
by Corollary~\ref{cor:UL}, $\UU=\LL$ is a global attractor in
$\normLone{\cdot}$.
\end{proof}
\begin{rem}
The condition given in Theorem \ref{thm:EndPointCriterion} is not
necessary. When $\mu$ is the Lebesgue measure it is shown in
Remark~\ref{rem:LowL}.\qed
\end{rem}

\section{Examples}\label{sec:Examples}
We consider three examples of the measure $\mu$ which defines the
convolution $\Phi$: the uniform distribution $d\mu(x)=dx$ on the
interval $[0,1]$, the exponential distribution
$d\mu(x)=ae^{-ax}dx$ on the segment $\Real^{+}$ and the Dirac
measure at a point $t\in COInt{0,1}$, \emph{i.e.,}
$\mu(A)=\delta_t(A)=\chi(A)(t)$, where again $\chi$ is an
indicator function of the set $A$ at the point  $t$.
\begin{thm}\label{thm:examples}
\ \\
\textbf{Uniform distribution on $[0,1]$}. For every $\lambda>0$,
the map $\TT$ has a unique fixed point:
 \[
  \breve{F}_A(x) =
\displaystyle\frac{A e^{A \lambda x}}{e^{A \lambda} ( A-1) + e^{A
\lambda x}}\,,
  \]
 which is a global attractor. Here $1\le A = A(\lambda)$ is
 the unique solution of the equation:
 \[
 e^{A\lambda}(A-1)^2=1\,.
 \]
\textbf{Exponential distribution}
(\emph{cf.}\cite{gamarnikMaxWeightIndSet}).
    For every $\lambda>0$ and every $a$,  the map $\TT$ has a unique
fixed point:
\[
\check{F}_K(x)=\exp(-Ke^{-ax})\,,
\]
which is a global attractor. Here $0\le K=K(\lambda)$ (independent
on $a$) is the unique solution of $\check{f}(K)=K$:
\[
\check{f}(K)=\lambda\frac{1-e^{-K}}{K}\,.
\]
\textbf{Dirac measure $\delta_t$}(\emph{cf.} \cite{KarpSipser}).
 For every $0<\lambda\le e$ the map $\TT$
has a unique fixed point:
\[
\tilde{F}_M(x)=\left\{%
\begin{array}{ll}
    0, &\quad x<0  \\
    M, &\quad 0<x<t \\
    1, & \quad t<x
\end{array}%
\right.
\]
which is a global attractor. Here $M$ is a unique solution of
$\tilde{f}(M)=M$:
\[
\tilde{f}(M)=e^{-\lambda M}\,.
\]
For every $\lambda>e$ there is no fixed point which is a global
attractor.
\end{thm}
\begin{proof}
In each case we first investigate the properties of the equations
on the parameters.
\subsection*{The Lebesgue measure}\label{sub:Leb} In this
case:
\[
\Phi(F)(x)=\int_0^{1-x}F(z) dz.
\]
\begin{prop}\label{prop:LebPeriod2}
The only twice differentiable function satisfying $\TT^2(F) = F$,
$F(1)=1$ and $F$ non decreasing  is the function ${\breve{F}}$
with $A\ge 1$. It follows that $\breve{F}_A$ is a fixed point of
$\TT$.
\end{prop}
\begin{proof}
We want to solve the equation $\TT^2(F) = F$, i.e.
\[
e^{\displaystyle {-\lambda \int_{0}^{1-x} e^{-\lambda
\int_{0}^{1-z} F(s) ds} dz }} = F(x)\,.
\]
After applying $\ln$ to both sides of the above equation and
differentiating them with respect to $x$ we get
 \beq \label{eqnLebDiffeq2}
  \lambda e^{-\lambda \int_{0}^{x} F(s) \,ds} =
  \displaystyle\frac{F'(x)}{F(x)}\,.
 \eeq
We repeat the same procedure once again, and we obtain the
second-order differential equation:
\[
- \lambda F(x) =
\displaystyle\frac{F''(x)F(x)-(F'(x))^2}{F'(x)F(x)}\,,
\]
in other words
\[
F'' F - (F')^2+\lambda F' F^2  = 0\,.
\]
This equation does not contain the independent variable, so the
standard substitution $z = F'(F)$ allows us to lower the degree of
the equation to $1$. Easily we get the solution $F(x) = \frac{C A
e^{A \lambda x}}{1 + C e^{A \lambda x}}$.  From $F(1) = 1$ there
follows $C = \frac{e^{-A \lambda}}{A-1}$, and
\[
F(x) = \displaystyle\frac{A e^{A \lambda x}}{e^{A \lambda} ( A-1)
+ e^{A \lambda x}}\,,
\]
and from $F(x)\ge0$, $F'(x)\ge 0$ on $[0,1]$ there follows $A\ge
1$. We substitute this function into the equation
\eqref{eqnLebDiffeq2}, and obtain:
\[
\displaystyle\frac{(1+(A-1) e^{A \lambda}) \lambda}
  {(A-1) e^{A \lambda} + e^{A \lambda x}}
=
  \frac{(A-1) A e^{A \lambda} \lambda}{(A-1) e^{A \lambda} + e^{A \lambda x}},
  \]
so $(A-1)^2 e^{A \lambda} = 1$. By assumption $F=\breve{F}_A$ is a
periodic point of period two. Because $\TT$ preserves both
additional conditions $\TT\breve{F}_A$ fulfils the assumptions of
the proposition and hence $\TT\breve{F}_A=\breve{F}_A$.
\end{proof}
The  monotone maps $\LL$ and $\UU$ are the images of themselves
under the map $\TT$ which is a composition of a smooth exponential
map with the  convolution $\Phi$ with the  smooth kernel $1$. That
means that $\LL$ and $\UU$ are at least twice differentiable and
hence satisfy the assumption of Proposition~\ref{prop:LebPeriod2},
so must be both equal to the map~$\breve{F}$.
\begin{rem}\label{rem:LebEqFixPoint} It is possible to construct and solve a
differential equation for the fixed point of $\TT$. However as the
convolution $\Phi$ exchanges the argument $x$ into $1-x$ one gets:
\[
-\lambda \frac{d}{dx}\Phi F (x) = \lambda F(1-x) =
\frac{F'(x)}{F(x)}\,,
\]
which is not a differential equation. One has to use special
symmetries to get rid of $F(1-x)$. Indeed, integrating and using
$F(1) = 1$ one obtains:
\[
F(x) = \lambda \int_0^x F(z) F(1-z) dz + 1-k\,,
\]
which produces the needed formula for $F(1-x)$.\qed
\end{rem}
\begin{rem}\label{rem:LowL}
After simplifications we have $\breve{F}_A(0)=A-1$, where $A$ runs
from one to two as $\lambda$ runs from infinity to zero. For
$\lambda>e/(e+1)$ we have $A-1<1/e$ which shows that the condition
in Theorem~\ref{thm:EndPointCriterion} is not necessary. \qed
\end{rem}
\subsection*{The exponential measure}\label{sub:Exp}
Here we work with a slightly different setting, namely the
interval of the arguments  in the definition of $\Dol$ becomes the
segment $(0,+\infty)$. There is a simplification in the
convolution:
\[
\Phi(F)(x)=\int_x^\infty F(z-x)a e^{-az}\,dz=e^{-ax}\int_0^\infty
F(w)ae^{-aw}\,dw=e^{-ax}\E_a[F]\,,
\]
which means that dependance on $x$ is outside the integral, so
that we can write:
 \[
 \TT(F)(x)=\exp(-\lambda e^{-ax}
 \int_{0}^{\infty}F(z)ae^{-az}\,dz)=
 \exp(-\lambda e^{-ax}\E_a[F])\,.
 \]
 If we set $G(x)=F(x/a)$ then we get a conjugated evolution which is
 independent on $a$:
\[
 \check{\TT}(G)(x)=\exp(-\lambda
 e^{-x}\E[G]),\quad\text{where}\quad
\E[G]=\E_1[G]=\int_{0}^{\infty}G(z)e^{-z}\,dz \,.
 \]
 In other words after one iteration the collection of distributions
 consists of one parameter family of double exponential functions
 $\check{\Dol}=\set{\exp(-Ke^{-x}),K>0}$,
 and the dynamic in the space of functions
is reduced to the dynamics of the parameter $K$. For
$G\in\check\Dol$ with parameter $K$ we have $\check\TT
G\in\check\Dol$ with parameter $\lambda\E G$, or $K\mapsto \check
f(K)=\lambda\E G$:
\[
\check{f}(K)=\lambda\int_{0}^{\infty}G(z)e^{-z}\,dz=\lambda\int_0^1\exp(-Kt)\,dt
=\lambda\frac{1-e^{-K}}{K}\,.
\]
In particular the equation $ \check{f}(\check{f}(K))=K$ has a
unique positive solution (which is also a solution of
$\check{f}(K)=K$), which means that $\TT$ has no periodic points
of proper period two, or that $\LL=\UU$.
\begin{rem}
Also here for $\lambda>e/(e-1)$ we have $K>1$ and hence
$\check{F}(0)=e^{-K}<1/e$.
\end{rem}
\subsection{The Dirac measure}\label{sub:Dirac}
\begin{lem}
\label{lem:OneDim}
 For each $\lambda$ the real function $\tilde{f}(x)=e^{-\lambda x}$ has a unique fixed point $M =
 M(\lambda)$. For $\lambda \le e$, the point $M$ is an attractor, for $\lambda >
 e$ it is a repeller, but then $\tilde{f}$ has an attracting periodic orbit
which attracts every $x\not=M$.
\end{lem}
\begin{proof}
The function $f_\lambda (x)$ is decreasing, so the first statement
is obvious. Using Implicit Function Theorem one checks that
$dM/d\lambda<0$. We have $\tilde{f}'(M) <0$, and $(d/d \lambda)
(\tilde{f}'(M))<0$,
 so the derivative of $\tilde{f}$ at $M$ is also a decreasing function of $\lambda$.
The equation $\tilde{f}'(M) = -1$, that is $-\lambda e^{-\lambda
M} = -\lambda M = -1$ has a unique solution with $\lambda=e$ and
$M = 1/e$. In order to check that for $\lambda\le e$ this  gives
the global attractor consider the second iterate, and see that for
$\lambda\le e$, $(\tilde{f}(\tilde{f}(x))-x)(M-x)>0$ for
$x\not=M$, and that for $\lambda>e$ there is an orbit of period 2
attracting all $x\not=M$. The details of this exercise are
omitted.
\end{proof}
The explicit form of $\TT F$ in the case of $\mu=\delta_b$ is
given by:
\[
\TT F(x) = \exp(-\lambda \int_x^1 F(z-x) \delta_b(z)) = \left\{
\begin{array}{ll}
e^{-\lambda F(b-x)} &  x \in[0,b]\\
1 & x \in \OCInt{b,1}
 \end{array} \right.\,.
\]
We observe that if $F$ is constant on the interval $[0,b]$ so is
$\TT(F)$. One can easily check that
\[
\One^n(x) = \left\{
\begin{array}{ll}
\tilde{f}^{n+1}(0) &  x \in[0,b]\\
1 & x \in\OCInt{ b,1}
 \end{array} \right. ,
\]
and the theorem follows from Lemma~\ref{lem:OneDim}.
\begin{rem}\label{rem:DeltaOtherAttractors}
For  $\lambda>e$  the points $L<M<U$ are the two periodic points
of period two which attracts (almost) every trajectory of
$\tilde{f}$. For $F\le M$ we have $\TT^{2n} F\to L$ and
$\TT^{2n+1} F\to U$ (for $F\ge U$ we have an analogous statement)
point-wise at every point $x\le b$ . If $F(x)=M$ for some $x\le b$
then $\TT^{2n}F(x)=M$ and $\TT^{2n+1}F(b-x)=M$.  Let $J$ be a
minimal interval which contains all such points $x$ and $a-x$,
then outside $J$ the iterates accumulate on $L$ and $U$ while
inside $J$ the accumulation points swap between $L$, $M$ and $U$
depending on the values of $F$  at points $a-x$ (but preserving
the monotony of $\TT^nF$). This shows that there is no simple
attracting point for all $F\in\Dol$ and shows that the bound in
Theorem~\ref{thm:TheECut} is tight. The details are skipped.\qed
\end{rem}
\begin{rem} The dynamics get much more complicated already for
$\mu=p\delta_{a}+(1-p)\delta_b$. There is a partition of $[0,1]$
into $N=2/(b-a)$ intervals of the points $m(b-a)$ and $b-m(b-a)$,
$m=0,1\dots$ such that a subset of $\Dol$ of maps constant on
these intervals is invariant under $\TT$ (independently on
$\lambda$). Studying this invariant subset is sufficient as that
is where all the iterates of $\One$ live. This reduces the
dynamics in the functional space into the dynamics in  $\Real^N$.
Similarly for $\mu=\sum p_i\delta_i$. In the simplest case
$a_i=(i-1)/n$ the invariant partition consists of $n$ intervals,
and the resulting dynamical system is a map from $[0,1]^n$ into
itself.
 \qed
\end{rem}
This concludes the description of the examples and the proof of
Theorem~\ref{thm:examples}.
\end{proof}

\bibliographystyle{amsalpha}

\providecommand{\bysame}{\leavevmode\hbox
to3em{\hrulefill}\thinspace}
\providecommand{\MR}{\relax\ifhmode\unskip\space\fi MR }
\providecommand{\MRhref}[2]{%
  \href{http://www.ams.org/mathscinet-getitem?mr=#1}{#2}
} \providecommand{\href}[2]{#2}

\end{document}